\documentclass{amsart}
\usepackage{amsfonts}
\usepackage{amsmath}
\usepackage{amssymb}
\usepackage{amsthm}

\newtheorem{theorem}{Theorem}

\newtheorem{lemma}{Lemma}

\begin{document}

\title[Apartments preserving transformations of Grassmannians]{Apartments preserving transformations of Grassmannians of infinite-dimensional vector spaces}
\author{Mark Pankov}
\address{Department of Mathematics and Computer Science, 
University of Warmia and Mazury, S{\l}oneczna 54, Olsztyn, Poland}
\email{pankov@matman.uwm.edu.pl}

\maketitle

\begin{abstract}
 We define the Grassmannians of an infinite-dimensional vector space $V$ 
 as the orbits of the action of the general linear group ${\rm GL}(V)$ on the set of all subspaces.
Let ${\mathcal G}$ be one of these Grassmannians. 
An apartment in ${\mathcal G}$ is the set of all elements of ${\mathcal G}$ spanned by subsets of a certain basis of $V$.
We show that every bijective transformation $f$ of ${\mathcal G}$ such that $f$ and $f^{-1}$ send apartments to apartments 
is induced by a semilinear automorphism of $V$.
In the case when ${\mathcal G}$ consists of subspaces whose dimension and codimension both are infinite,
a such kind result will be proved also for the connected components of the associated Grassmann graph.
\end{abstract}

\section{Introduction}
The first geometric characterization of semilinear  isomorphisms was the Fundamental Theorem of Projective Geometry 
which states that all isomorphisms between projective spaces are induced by semilinear isomorphisms of the associated vector spaces. 
Note that this statement holds for vector spaces of an arbitrary (not necessarily finite) dimension.
The second breakthrough in this direction was Chow's theorem  describing automorphisms of 
the Grassmann graphs of finite-dimensional vector spaces \cite{Chow}, see also \cite{Pankov-book1,Wan}.

Other characterization of semilinear isomorphisms is related to the concept of {\it apartment}.  
Every apartment in a Grassmannian of a finite-dimensional vector space 
consists of all elements of this Grassmannian spanned by subsets of a certain basis of the vector space,
in other words, this is the intersection of the Grassmannian with an apartment in the Tits building  for the associated general linear group \cite{Tits}.
By \cite{Pankov-book1}, all apartments preserving bijections between Grassmannians of finite-dimensional vector spaces are induced by semilinear isomorphisms.
In the present paper, the same statement  will be established for Grassmannians of infinite-dimensional vector spaces.

Grassmannians of an infinite-dimensional vector space can be defined as the orbits of the action of the general linear group on the set of all subspaces. 
There are the following two types of Grassmannians:
\begin{enumerate}
\item[$\bullet$] Grassmannians formed by subspaces of finite dimension or codimension,
\item[$\bullet$] Grassmannians consisting of subspaces whose dimension and codimension both  are infinite,
in particular, the Grassmannian formed by subspaces whose dimension and codimension are equal to the dimension of the vector space.
\end{enumerate}
Grassmannians of first type are similar to Grassmannians of finite-dimensional vector spaces.
Grassmannians of second type are more complicated: 
they contain infinite subsets of mutually incident elements and the associated Grassmann graphs are not connected.
Apartments are defined as for Grassmannians of finite-dimensional vector spaces, i.e. by bases of vector spaces.

Consider the Grassmannian consisting of subspaces whose dimension and codimension are equal to the dimension of the vector space.
This Grassmannian contains pairs of complementary subspaces. 
A. Blunck and H. Havlicek \cite{BH} proved that 
every bijective transformation preserving the complementarity of subspaces is an automorphism of the corresponding Grassmann graph.

All automorphisms of the Grassmann graphs related to Grassmannians of first type 
are induced by semilinear automorphisms of the corresponding vector space or the dual vector space.
For the Grassmann graphs associated to Grassmannians of second type this fails.
By L. Plevnik \cite{Plevnik},  
the restrictions of automorphisms of Grassmann graphs to connected components can be obtained from pairs of semilinear isomorphisms of special type
and there are graph automorphisms whose restrictions to connected components are not induced by semulinear isomorphisms.
Also, using  Blunck--Havlicek's result mentioned above,  Plevnik showed that every complementarity preserving bijection is induced by a semilinear isomorphism. 

Our first result states that all apartments preserving transformations of Grassmannians of an infinite-dimensional vector space are induced by 
semilinear automorphisms of this vector spaces.
In \cite{Pankov}, this statement was proved for Grassmannians formed by subspaces whose dimension is less than the dimension of the vector space.
As in \cite{Pankov}, we will use local methods related to special bijections between apartments.
The second result concerns apartments preserving bijective transformations of connected components of Grassmann gtaphs.

\section{Results} 
Let $V$ be a left vector space over a division ring. 
Throughout the paper we suppose that the dimension of $V$ is an infinite cardinal number $\alpha$.
For every cardinal number $\beta \le \alpha$ we denote by ${\mathcal G}_{\beta}(V)$ 
the Grassmannian consisting of all subspaces of $V$ whose dimension is $\beta$ and whose codimension is $\alpha$.
Also, we write ${\mathcal G}^{\beta}(V)$ for the Grassmannian formed by all subspaces of $V$ whose dimension is $\alpha$ and whose codimension is $\beta$.
Then ${\mathcal G}_{\alpha}(V)={\mathcal G}^{\alpha}(V)$ consists of all subspaces whose dimension and codimension is $\alpha$.
It was noted above that each Grassmannian is an orbit of the action of the group ${\rm GL}(V)$ on the set of all subspaces of $V$.

Let ${\mathcal G}$ be a Grassmannian of $V$, i.e. 
${\mathcal G}$ coincides with ${\mathcal G}_{\beta}(V)$ or ${\mathcal G}^{\beta}(V)$ for a certain cardinal number $\beta\le \alpha$.
For every basis $B$ of the vector space $V$ the set of all elements of ${\mathcal G}$ spanned by subsets of $B$ is called 
the {\it apartment} of ${\mathcal G}$ defined by the basis $B$.
If two bases of $V$ define the same apartment, then the vectors from one of the bases are scalar multiples of the vectors from the other.
It is not difficult to prove that for any two elements of ${\mathcal G}$ there is an apartment containing them
\cite[Proposition 1.4]{Pankov-book1}.

A subset of ${\mathcal G}_{1}(V)$ is an apartment if and only if  it is a basis of the projective space $\Pi_{V}$ associated to the vector space $V$. 
Every apartment of ${\mathcal G}^{1}(V)$ is a proper subset in a basis of the dual projective space $\Pi^{*}_{V}$.

A bijective transformation $l:V\to V$ is a {\it semilinear automorphism} if 
$$l(x+y)=l(x)+l(y)$$ 
for all vectors $x,y\in V$ and there is an automorphism $\sigma$ of the division ring associated to $V$ such that 
$$l(ax)=\sigma(a)l(x)$$
for every vector $x\in V$ and every scalar $a$.
Each semilinear automorphism of $V$ induces a bijective transformation of the Grassmannian ${\mathcal G}$
preserving the family of apartments. 
Two semilinear automorphisms define the same transformation of ${\mathcal G}$ if and only if one of them is a scalar multiple of the other.

\begin{theorem}\label{theorem1}
If $f$ is a bijective transformation of ${\mathcal G}$ such that $f$ and $f^{-1}$ send apartments to apartments,
then $f$ is induced by a semilinear automorphism of $V$.
\end{theorem}

In the case when ${\mathcal G}={\mathcal G}_{1}(V)$, this statement is a simple consequence of the Fundamental Theorem of Projective Geometry.
Indeed, three distinct points of the projective space $\Pi_{V}$ are non-collinear if and only if there is an apartment of ${\mathcal G}_{1}(V)$ containing them. 
This means that $f$ and $f^{-1}$ transfer triples of non-collinear points to triples of non-collinear points.
The latter is equivalent to the fact that    $f$ and $f^{-1}$ send triples of collinear points to triples of collinear points.
Therefore, $f$ is an automorphism of $\Pi_{V}$ and the Fundamental Theorem of Projective Geometry implies that 
it is induced by a semilinear automorphism of $V$.

Elements $X,Y\in {\mathcal G}$ are called {\it adjacent} if 
$$\dim X/(X\cap Y)=\dim Y/(X\cap Y)=1,$$
in other words, $X\cap Y$ is a hyperplane in both $X,Y$. 
The {\it Grassmann graph} $\Gamma({\mathcal G})$ is the graph whose vertex set is ${\mathcal G}$ and  whose edges are pairs of adjacent elements.
The graph is connected if $\beta$ is finite, and any two distinct vertices of this graph are adjacent if $\beta=1$.
In the case when $\beta$ is infinite, the Grassmann graph is not connected and the connected component containing $X\in {\mathcal G}$ consists of all $Y\in {\mathcal G}$ such that
\begin{equation}\label{eq1}
\dim X/(X\cap Y)=\dim Y/(X\cap Y)<\infty.
\end{equation}
The bijective transformations of ${\mathcal G}$ induced by semilinear automorphisms of $V$ are automorphisms of the Grassmann graph $\Gamma({\mathcal G})$.

If $\beta$ is finite and greater than $1$, then all automorphisms of the Grassmann graphs associated to ${\mathcal G}_{\beta}(V)$ and ${\mathcal G}^{\beta}(V)$
are  induced by semilinear automorphisms of $V$ and $V^{*}$, respectively.
The proof of this statement is based on the intersection properties of maximal cliques in Grassmann graphs and similar to the proof of classical Chow theorem.

Suppose that $\beta$ is infinite. 
The restrictions of an automorphism of $\Gamma({\mathcal G})$ to distinct connected components can be induced by distinct semilinear automorphisms of $V$
and there are automorphisms of $\Gamma({\mathcal G})$ whose restrictions to connected components are not induced by semilinear automorphisms of $V$  
\cite[Corollary 7]{Plevnik}.
By \cite[Theorem 1]{Plevnik}, the restriction of every automorphism of $\Gamma({\mathcal G})$ to any connected component 
can be obtained from a pair of semilinear isomorphisms of special type.

Let ${\mathcal C}$ be a connected component of $\Gamma({\mathcal G})$ (as above, we suppose that $\beta$ is infinite).
If ${\mathcal A}$ is an apartment of ${\mathcal G}$ having a non-empty intersection with ${\mathcal C}$ and $X\in {\mathcal A}\cap {\mathcal C}$,
then ${\mathcal A}\cap {\mathcal C}$ consists of all $Y\in {\mathcal A}$ satisfying \eqref{eq1}.
Every such intersection will be called an {\it apartment} of the component ${\mathcal C}$.

\begin{theorem}\label{theorem2}
Suppose that $\beta$ is infinite and ${\mathcal C}, {\mathcal C}'$ are connected components of $\Gamma({\mathcal G})$.
If $f$ is a bijection of ${\mathcal C}$ to ${\mathcal C}'$ such that $f$ and $f^{-1}$ send apartments to apartments, 
then $f$ is induced by a semiliear automorphism of $V$.
\end{theorem}

\section{Preliminary}

\subsection{Stars and tops}
For a subspace $X\subset V$ the set ${\mathcal G}^{1}(X)$ consisting of all hyperplanes of $X$ and  
the set ${\mathcal S}(X)$ consisting of all subspaces which contain $X$ as a hyperplane are called the {\it top} and the {\it star} associated to $X$, respectively.
Every star and every top is formed by mutually adjacent elements of a certain Grassmannian.
If $\beta>1$, then every maximal clique of the graph $\Gamma({\mathcal G})$ is a top or a star \cite[Sections 3.1 and 3.8]{Pankov-book1}.

Recall that every line in the projective space $\Pi_{V}$ consists of all $1$-dimensional subspaces contained in a certain element of ${\mathcal G}_{2}(V)$
and every line in the dual projective space $\Pi^{*}_{V}$ is formed by all hyperplanes of $V$ which contain a certain element of ${\mathcal G}^{2}(V)$.
Suppose that $\beta>1$ and ${\mathcal S}(X)$ and ${\mathcal G}^{1}(Y)$ are a star and a top in ${\mathcal G}$.
Their intersection consists of all $Z\in {\mathcal G}$ satisfying  $X\subset Z\subset Y$. The intersection is called a {\it line} of ${\mathcal G}$ if it is non-empty.
Every star ${\mathcal S}(X)$ (together with all lines contained in it) can be identified with the projective space $\Pi_{V/X}$ and
every top ${\mathcal G}^{1}(Y)$ (together with all lines contained in it) is the projective space $\Pi^{*}_{Y}$.
The intersection of two distinct maximal cliques is empty or a one-element set or a line.
Therefore, every automorphism of  $\Gamma({\mathcal G})$ sends lines to lines and
its restriction to a maximal clique (a star or a top) is an isomorphism of the corresponding projective space to the projective space associated to other maximal clique.
If an automorphism of $\Gamma({\mathcal G})$ transfers a star ${\mathcal S}(X)$ to a top ${\mathcal G}^{1}(Y)$, 
then the corresponding projective spaces are isomorphic and we have $\dim (V/X)=\dim Y^{*}$.
The latter equality implies that $\dim Y^{*}$ is infinite and, by \cite[Section II.3]{Baer}, it is equal to ${\gamma}^{\dim Y}$,
where $\gamma$ is the cardinal number of the associated division ring.
This means that every automorphism of $\Gamma({\mathcal G})$ sends stars to stars and tops to tops except the case when 
${\mathcal G}={\mathcal G}_{\beta}(V)$, $\beta$ is infinite and $\alpha={\gamma}^{\beta}$.

Suppose that $\beta$ is infinite and ${\mathcal C}$ is a connected component of the graph $\Gamma({\mathcal G})$.
Let $f$ be the restriction of an automorphism of $\Gamma({\mathcal G})$ to ${\mathcal C}$ which sends stars to stars and tops to tops.
It is clear that $f$ transfers ${\mathcal C}$ to a certain connected component ${\mathcal C}'$.
Denote by ${\mathcal C}_{-1}$ and ${\mathcal C}'_{-1}$ the sets formed by all hyperplanes in elements of ${\mathcal C}$
and all hyperplanes in elements of ${\mathcal C}'$, respectively.
Since ${\mathcal C}_{-1}$ and ${\mathcal C}'_{-1}$ can be identified with the sets of all stars in ${\mathcal C}$ and ${\mathcal C}'$ (respectively),
$f$ induces a bijection $f_{-}:{\mathcal C}_{-1}\to {\mathcal C}'_{-1}$.
Denote by ${\mathcal C}_{+1}$ and ${\mathcal C}'_{+1}$ 
the sets formed by all subspaces which contain as a hyperplane an element of ${\mathcal C}$ or an element of ${\mathcal C}'$, respectively.
These sets can be considered as the sets of all tops in ${\mathcal C}$ and ${\mathcal C}'$, respectively.
Therefore, $f$ induces a bijection $f_{+}:{\mathcal C}_{+1}\to {\mathcal C}'_{+1}$.
Observe that ${\mathcal C}_{-1},{\mathcal C}_{+1},{\mathcal C}'_{-1},{\mathcal C}'_{-1}$ are connected components of the graph $\Gamma({\mathcal G})$.
Two distinct stars in ${\mathcal C}$ or ${\mathcal C}'$ have a non-empty intersection if and only if 
the corresponding elements of ${\mathcal C}_{-1}$ or ${\mathcal C}'_{-1}$ are adjacent. The same holds for tops.
This means that $f_{-}$ and $f_{+}$ are adjacency preserving in both directions, i.e. these are isomorphisms between the restrictions of the graph $\Gamma({\mathcal G})$
to the corresponding connected components.
We apply the above arguments to $f_{-}$ and $f_{+}$. Step by step, we extend $f$ to a bijection ${\overline f}:{\mathcal C}_{\pm}\to {\mathcal C}'_{\pm}$, where
${\mathcal C}_{\pm}$ is the set of all $X\in {\mathcal G}$ such that $X$ is a subspace of finite codimension in an element of ${\mathcal C}$ or 
$X$ contains an  element of ${\mathcal C}$ as a subspace of finite codimension
and ${\mathcal C}'_{\pm}$ is similarly defined for ${\mathcal C}'$. 
Note that ${\overline f}$ is inclusion preserving in both directions, 
i.e. $X\subset Y$ if and only if ${\overline f}(X)\subset {\overline f}(Y)$ for all $X,Y\in {\mathcal C}_{\pm}$.
By \cite{Plevnik}, ${\overline f}$ is defined by a pair of special semilinear isomorphisms, see the case (A) in \cite[Theorem 1]{Plevnik} for the details. 

\subsection{Induced transformations}
Let ${\mathcal G}_{1}$ and ${\mathcal G}_{2}$ be Grassmannians of $V$.
Suppose that a bijective transformation $f_{1}$ of ${\mathcal G}_1$ {\it induces} a bijective transformation $f_{2}$ of ${\mathcal G}_2$, i.e.
$X\in {\mathcal G}_1$ is incident to $Y\in {\mathcal G}_2$ if and only if $f_{1}(X)$ is incident to $f_{2}(Y)$.

If $f_1$ is induced by a semilinear automorphism $l$ of $V$, i.e. $f_{1}(X)=l(X)$ for every $X\in {\mathcal G}_{1}$,
then $X\in {\mathcal G}_1$ is incident to $Y\in {\mathcal G}_2$ if and only if $l(X)$ is incident to $f_{2}(Y)$.
This implies that $f_{2}(Y)=l(Y)$ for every $Y\in {\mathcal G}_{2}$.

We say that apartments ${\mathcal A}_{1}\subset {\mathcal G}_1$ and  ${\mathcal A}_{2}\subset {\mathcal G}_2$
are {\it associated} if they are defined by the same basis of $V$.
If $f_1$ sends apartments to apartments, then $f_{2}$ transfers the apartment of ${\mathcal G}_{2}$ associated to 
an apartment ${\mathcal A}\subset {\mathcal G}_1$ to the apartment associated to $f_{1}({\mathcal A})$.

So, we get the following.

\begin{lemma}\label{lemma0-2}
Suppose that a bijective transformation of a Grassmannian induces a bijective transformation of other Grassmannian.
Then the following assertions are fulfilled:
\begin{enumerate}
\item[(1)] if one of these transformations is induced by a semilinear automorphism, 
then the other transformation is induced by the same semilinear automorphism;
\item[(2)] if one of these transformations sends apartments to apartments, then the same holds for the other.
\end{enumerate}
\end{lemma}

\section{Proof of Theorem \ref{theorem1}}
Recall that ${\mathcal G}$ coincides with ${\mathcal G}_{\beta}(V)$ or ${\mathcal G}^{\beta}(V)$, where $\beta$ is a cardinal number not greater than $\alpha$.
If $\beta=\alpha$, then ${\mathcal G}={\mathcal G}_{\alpha}(V)={\mathcal G}^{\alpha}(V)$.
Also, we suppose that $f$ is a bijective transformation of ${\mathcal G}$ such that $f$ and $f^{-1}$ send apartments to apartments.

\subsection{The case $\beta=1$}
It was noted above that the statement is a simple consequence of the Fundamental Theorem of Projective Geometry if ${\mathcal G}={\mathcal G}_{1}(V)$. 
Consider the case when ${\mathcal G}={\mathcal G}^{1}(V)$ and show that $f$ induces a bijective transformation of ${\mathcal G}_{1}(V)$, 
i.e. for every $1$-dimensional subspace $P\subset V$ there is a $1$-dimensional subspace $P'\subset V$ such that 
a hyperplane $X\subset V$ contains $P$ if and only if $f(X)$ contains $P'$.

For a $1$-dimensional subspace $P\subset V$ we take any apartment of ${\mathcal G}_{1}(V)$ containing it
and denote by ${\mathcal A}$ the associated apartment of ${\mathcal G}^{1}(V)$.
There is the unique element $X\in{\mathcal A}$ which does not contain $P$. 
Since $f({\mathcal A})$ is an apartment of ${\mathcal G}^{1}(V)$, 
the intersection of all elements from $f({\mathcal A}\setminus \{X\})$ is a $1$-dimensional subspace $P'$.
If a hyperplane $Y\subset V$ does not contain $P$, then 
$$({\mathcal A}\setminus \{X\})\cup \{Y\}$$ 
is an apartment of ${\mathcal G}^{1}(V)$ and the same holds for 
$$f(({\mathcal A}\setminus \{X\})\cup \{Y\})$$
which implies that $f(Y)$ does not contain $P'$. 
We apply the same arguments to $f^{-1}$ and establish that a hyperplane $Y\subset V$ does not contain $P$ if and only if $f(Y)$ does not contain $P'$.
This means that  $P'$ is as required.

So, $f$ induces a bijective transformation $g$ of ${\mathcal G}_{1}(V)$.
By the second part of Lemma \ref{lemma0-2}, $g$ and $g^{-1}$ send apartments to apartments.
Then $g$ is induced by a semilinear  automorphism of $V$.
It follows from the first part of Lemma \ref{lemma0-2} that $f$ is induced by the same semilinear automorphism.

\subsection{Special bijections between apartments}
From this moment, we suppose that $\beta>1$.
Let $\{e_i\}_{i\in I}$ be a basis of $V$ and let ${\mathcal A}$ be the associated apartment of ${\mathcal G}$.
For every $i\in I$ we denote by ${\mathcal A}(+i)$ the subset of ${\mathcal A}$ consisting of all elements which contain $e_{i}$
and write ${\mathcal A}(-i)$ for the subset of ${\mathcal A}$ formed by all elements which do not contain $e_{i}$.
We say that ${\mathcal A}(+i)$ and ${\mathcal A}(-i)$ are {\it simple subsets of first} and {\it second type}, respectively.
For distinct $i,j\in I$ we define
$${\mathcal A}(+i,+j)={\mathcal A}(+i)\cap {\mathcal A}(+j),$$
$${\mathcal A}(+i,-j)={\mathcal A}(+i)\cap {\mathcal A}(-j).$$
A subset ${\mathcal X}\subset{\mathcal A}$ is called {\it inexact}
if there is an apartment of ${\mathcal G}$ distinct from ${\mathcal A}$ and containing ${\mathcal X}$.
For any distinct $i,j\in I$ the subset 
\begin{equation}\label{eq2}
{\mathcal A}(+i,+j)\cup {\mathcal A}(-i)
\end{equation}
is inexact and every maximal inexact subset of ${\mathcal A}$ is a subset of type \eqref{eq2}, see \cite[Sections 3.4 and 3.8]{Pankov-book1}.

We say that ${\mathcal X}\subset {\mathcal A}$ is a {\it complementary} subset if 
${\mathcal A}\setminus {\mathcal X}$ is a maximal inexact subset.
The complementary subset corresponding to the maximal inexact subset \eqref{eq2} is ${\mathcal A}(+i,-j)$.
For two distinct complementary subsets ${\mathcal A}(+i,-j)$ and ${\mathcal A}(+i',-j')$ one of the following possibilities is realized:
\begin{enumerate}
\item[$\bullet$] $i=i'$ or $j=j'$,
\item[$\bullet$] $i=j'$ or $j=i'$, i.e. the intersection of the complementary subsets is empty,
\item[$\bullet$] $\{i,j\}\cap \{i',j'\}=\emptyset$.
\end{enumerate}
In the first case, the complementary subsets are said to be {\it adjacent}.
Observe that distinct complementary subsets ${\mathcal X},{\mathcal Y}\subset {\mathcal A}$ are adjacent if and only if their intersection is maximal,
i.e. for any distinct complementary subsets ${\mathcal X}',{\mathcal Y}'\subset {\mathcal A}$
satisfying
$${\mathcal X}\cap{\mathcal Y}\subset {\mathcal X}'\cap{\mathcal Y}'$$
the inverse inclusion holds.
An easy verification shows that every maximal collection of mutually adjacent complementary subsets of ${\mathcal A}$ is
$$\{{\mathcal A}(+i,-j)\}_{j\in\, I\setminus \{i\}}\;\mbox{ or }\;\{{\mathcal A}(+j,-i)\}_{j\in\, I\setminus \{i\}}$$
for a certain $i\in I$.
This means that every simple subset of ${\mathcal A}$ can be characterized as the union of all elements from 
a maximal collection of mutually adjacent complementary subsets.

Let ${\mathcal A}'$ be the apartment of ${\mathcal G}$ defined by a basis $\{e'_i\}_{i\in I}$.
The simple subsets of ${\mathcal A}'$ associated to the vector $e'_{i}$ will be denoted by ${\mathcal A}'(+i)$ and ${\mathcal A}'(-i)$.

A bijection $g:{\mathcal A}\to {\mathcal A}'$ is called {\it special} if $g$ and $g^{-1}$ transfer inexact subsets to inexact subsets.

\begin{lemma}\label{lemma1-1}
If $g:{\mathcal A}\to {\mathcal A}'$ is a special bijection, 
then there exists a bijective transformation $\delta:I\to I$ such that one of the following possibilities is realized:
\begin{enumerate}
\item[(1)]
$g({\mathcal A}(+i))={\mathcal A}'(+\delta(i))$ and $g({\mathcal A}(-i))={\mathcal A}'(-\delta(i))$ for all $i\in I$,
\item[(2)]
$g({\mathcal A}(+i))={\mathcal A}'(-\delta(i))$ and $g({\mathcal A}(-i))={\mathcal A}'(+\delta(i))$ for all $i\in I$.
\end{enumerate}
\end{lemma}

This lemma is taken from \cite{Pankov}, see also \cite[Lemma 3.11]{Pankov-book1}.
We recall the proof, since the same arguments will be exploited to prove Theorem \ref{theorem2}. 

\begin{proof}
It is clear that $g$ and $g^{-1}$ send  maximal inexact subsets to maximal inexact subsets.
Then ${\mathcal X}\subset {\mathcal A}$ is a complementary subset if and only if $g({\mathcal X})$ is a complementary subset of ${\mathcal A}'$.
Since two complementary subsets are adjacent only in the case when their intersection is maximal, 
$g$ and $g^{-1}$ transfer adjacent complementary subsets to adjacent complementary subsets.
Every simple subset can be characterized as the union of all elements from  a maximal collection of mutually adjacent complementary subsets.
This guarantees that $g$ and $g^{-1}$ send simple subsets to simple subsets.
Two distinct simple subsets are of different types if and only if their intersection is empty or a complementary subset.
Therefore, $g$ and $g^{-1}$ transfer simple subsets of different types to simple subsets of different types;
in other words, they preserve the types of all simple subsets or change the type of every simple subset.
Since ${\mathcal A}(-i)={\mathcal A}\setminus {\mathcal A}(+i)$, for every $i\in I$ there exists $\delta(i)\in I$ such that
the case (1) or (2) is realized.
\end{proof}

We say that $g$ is a special bijection of {\it first} or {\it second type} if it satisfies (1) or (2), respectively.
Suppose that $X\in {\mathcal A}$ is spanned by all vectors $e_{i}$ with $i\in J\subset I$. 
Then (1) implies that $g(X)$ is spanned by all $e'_{i}$ with $i\in \delta(J)$.
In the case (2), $g(X)$ is spanned by all $e'_{i}$ with $i\in I\setminus\delta(J)$.
Therefore, the second possibility can be realized only for $\beta=\alpha$, i.e. if ${\mathcal G}={\mathcal G}_{\alpha}(V)={\mathcal G}^{\alpha}(V)$.

\begin{lemma}\label{lemma1-2}
If $g:{\mathcal A}\to {\mathcal A}'$ is a special bijection, then the following assertions are fulfilled: 
\begin{enumerate}
\item[$\bullet$] $g$ is adjacency preserving in both directions, 
i.e. it is an isomorphism between the restrictions of the Grassmann graph $\Gamma({\mathcal G})$ to ${\mathcal A}$ and ${\mathcal A}'$.
\item[$\bullet$] If $\beta$ is infinite and $g$ is a special bijection of first type, then $g$ is inclusion preserving in both directions.
\end{enumerate}
\end{lemma}

\begin{proof}
This is a simple consequence of Lemma \ref{lemma1-1}.
\end{proof}

\begin{lemma}\label{lemma1-3}
The transformation $f$ is an automorphism of the graph $\Gamma({\mathcal G})$.
\end{lemma}

\begin{proof}
For any two elements of ${\mathcal G}$ there is an apartment containing them.
The restriction of $f$ to every apartment ${\mathcal A}\subset {\mathcal G}$ is a special bijection of ${\mathcal A}$ to $f({\mathcal A})$.
Therefore, the statement follows from the first part of Lemma \ref{lemma1-2}.
\end{proof}

\begin{lemma}\label{lemma1-4}
The restriction of $f$ to every apartment is a special bijection of first type. 
\end{lemma}

\begin{proof}
The statement is trivial if ${\mathcal G}$ coincides with ${\mathcal G}_{\beta}(V)$ or ${\mathcal G}^{\beta}(V)$ and $\beta<\alpha$.
Suppose that ${\mathcal G}={\mathcal G}_{\alpha}(V)={\mathcal G}^{\alpha}(V)$.
Lemma \ref{lemma1-2} states that $f$ is an automorphism of $\Gamma({\mathcal G})$ and, by Subsection 3.1, it transfers stars to stars and tops to tops. 
Let ${\mathcal A}$ be an apartment of ${\mathcal G}$.
We take any $X\in {\mathcal A}$ and consider the star ${\mathcal S}(X)$.
The transformation  $f$ sends ${\mathcal A}\cap {\mathcal S}(X)$  to the intersection of the apartment $f({\mathcal A})$ with a certain star. 
If the restriction of $f$ to ${\mathcal A}$ is a special bijection of second type, 
then $f({\mathcal A}\cap {\mathcal S}(X))$ is the intersection of $f({\mathcal A})$ with a top and we get a contradiction.
\end{proof}

\begin{lemma}\label{lemma1-5}
If $\beta$ is infinite, then $f$ is inclusion preserving in both directions.
\end{lemma}

\begin{proof}
For any two elements of ${\mathcal G}$ there is an apartment containing them.
By Lemma \ref{lemma1-4}, the restriction of $f$ to any apartment is a special bijection of first type.
The second part of Lemma \ref{lemma1-2} gives the claim.
\end{proof}

\subsection{The case when $\beta$ is finite and greater than $1$}
Suppose that ${\mathcal G}$ coincides with ${\mathcal G}_{\beta}(V)$ or ${\mathcal G}^{\beta}(V)$, where $\beta$ is finite and greater than $1$.
By Lemma \ref{lemma1-3}, $f$ is an automorphism of the Grassmann graph $\Gamma({\mathcal G})$.
This implies that $f$ is induced by a semilinear automorphism of $V$ if ${\mathcal G}={\mathcal G}_{\beta}(V)$.

Consider the case when ${\mathcal G}={\mathcal G}^{\beta}(V)$. 
By Subsection 3.1, $f$ and $f^{-1}$ transfer tops to tops, i.e.
there is a bijective transformation $f_{\beta-1}$ of ${\mathcal G}^{\beta-1}(V)$ such that 
$$f({\mathcal G}^{1}(X))={\mathcal G}^{1}(f_{\beta-1}(X))$$
for every $X\in {\mathcal G}^{\beta -1}(V)$.
Since two distinct tops of ${\mathcal G}^{\beta}(V)$ have a non-empty intersection if and only if the corresponding elements of ${\mathcal G}^{\beta-1}(V)$ are adjacent,
$f_{\beta-1}$ is an automorphism of the Grassmann graph associated to ${\mathcal G}^{\beta-1}(V)$.
If $\beta\ge 3$, then we apply the same arguments to $f_{\beta-1}$ and get a bijective transformation $f_{\beta-2}$ of ${\mathcal G}^{\beta-2}(V)$.
Step by step, we come to a bijective transformation $f_1$ of ${\mathcal G}^{1}(V)$ induced by $f$.
The second part of Lemma \ref{lemma0-2} implies that  $f_{1}$ and the inverse transformation send apartments to apartments.
Then $f_1$ is induced by a semilinear  automorphism of $V$ (see Subsection 4.1) and 
$f$ is induced by the same semilinear automorphism (the first part of Lemma \ref{lemma0-2}).

\subsection{The case when $\beta$ is infinite}
Suppose that ${\mathcal G}$ coincides with ${\mathcal G}_{\beta}(V)$ or ${\mathcal G}^{\beta}(V)$ and $\beta$ is an infinite cardinal number.

Let $X\in {\mathcal G}$.
By Lemma \ref{lemma1-5}, $f$ is inclusion preserving in both directions.
Therefore, $Y\in {\mathcal G}$ is a hyperplane of $X$ if and only if $f(Y)$ is a hyperplane of $f(X)$.
Denote by $f_{X}$ the restriction of $f$ to ${\mathcal G}^{1}(X)$. This is a bijection to ${\mathcal G}^{1}(f(X))$.
Every apartment ${\mathcal A}\subset {\mathcal G}^{1}(X)$ can be extended to an apartment ${\mathcal A}'\subset {\mathcal G}$.
Then $f_{X}({\mathcal A})$ is the intersection of the apartment $f({\mathcal A}')\subset {\mathcal G}$ with ${\mathcal G}^{1}(f(X))$.
This intersection is an apartment of ${\mathcal G}^{1}(f(X))$.
So, $f_{X}$ is a bijection of ${\mathcal G}^{1}(X)$ to ${\mathcal G}^{1}(f(X))$ transferring apartments to apartments.
The same arguments show that the inverse bijection of ${\mathcal G}^{1}(f(X))$ to ${\mathcal G}^{1}(X)$ sends apartments to apartments.
As in Subsection 4.1, we establish that $f_{X}$ is induced by a semilinear isomorphism $$l_{X}:X\to f(X).$$
This semilinear isomorphism is unique up to a scalar multiple.

Now, we take any $Y\in {\mathcal G}$ contained in $X$.
Observe that $Y$ is contained in $Z\in {\mathcal G}^{1}(X)$ if and only if $f(Y)$ is contained in $f(Z)=l_{X}(Z)$.
This guarantees that $f(Y)=l_{X}(Y)$.
By the same arguments, for every $Z\in {\mathcal G}^{1}(Y)$ we have $f_{Y}(Z)=f(Z)=l_{X}(Z)$.
This means that the semilinear isomorphism $l_{Y}$ coincides with the restriction of $l_{X}$ to $Y$ up to a scalar multiple.

For every $X\in {\mathcal G}$ we denote by $h_{X}$ the bijection of ${\mathcal G}_{1}(X)$ to ${\mathcal G}_{1}(f(X))$ induced by $l_{X}$
(this is an isomorphism between the projective spaces $\Pi_{X}$  and $\Pi_{f(X)}$).
If $Y\in {\mathcal G}$ is contained in $X$, then  the restriction of $h_{X}$ to ${\mathcal G}_{1}(Y)$ coincides with $h_{Y}$.

\begin{lemma}\label{lemma2-1}
If $X,Y\in {\mathcal G}$ have a non-zero intersection, then $h_{X}(P)=h_{Y}(P)$ for every $1$-dimensional subspace $P\subset X\cap Y$.
\end{lemma}

\begin{proof}
If $X\cap Y$ is an element of ${\mathcal G}$, then $$h_{X}(P)=h_{X\cap Y}(P)=h_{Y}(P).$$
In the case when $X\cap Y$  does not belong to ${\mathcal G}$, 
there exists $Z\in {\mathcal G}$ containing $X\cap Y$ and intersecting both $X,Y$ in some elements of ${\mathcal G}$.
Then $h_{X}(P)=h_{Z}(P)=h_{Y}(P)$.
\end{proof}

By Lemma \ref{lemma2-1}, 
there is a bijective transformation  $h$ of ${\mathcal G}_{1}(V)$ whose restriction to every ${\mathcal G}_{1}(X)$, $X\in {\mathcal G}$ coincides with $h_{X}$. 
This is an automorphism of the projective space $\Pi_{V}$, i.e. it is induced by a semilinear automorphism of $V$.
This semilinear automorphism  induces $f$.

\section{Proof of Theorem \ref{theorem2}}
Suppose that ${\mathcal G}$ coincides with ${\mathcal G}_{\beta}(V)$ or ${\mathcal G}^{\beta}(V)$ and $\beta$ is an infinite cardinal number.
Let ${\mathcal C}$ and ${\mathcal C}'$ be connected components of the Grassmann graph $\Gamma({\mathcal G})$.
Consider a bijection $f:{\mathcal C}\to {\mathcal C}'$ such that $f$ and $f^{-1}$ send apartments to apartments.
Recall that apartments in a connected component of $\Gamma({\mathcal G})$ are the non-empty intersections of this component and
apartments of ${\mathcal G}$. 
For any two elements of a connected component of $\Gamma({\mathcal G})$ there is an apartment of this component containing them
(this follows immediately from the same property for apartments in Grassmannians).

Let ${\mathcal A}_{c}$ be an apartment of ${\mathcal C}$. 
Then ${\mathcal A}_{c}={\mathcal A}\cap {\mathcal C}$, where ${\mathcal A}$ is an apartment of ${\mathcal G}$. 
A subset ${\mathcal X}\subset{\mathcal A}_{c}$ is called {\it inexact}
if there is an apartment of ${\mathcal C}$ distinct from ${\mathcal A}_{c}$ and containing ${\mathcal X}$.
A subset of ${\mathcal A}_{c}$ is inexact if and only if it is the intersection of ${\mathcal C}$ with an inexact subset of ${\mathcal A}$;
moreover, there is a one-to-one correspondence between 
maximal inexact subsets of ${\mathcal A}_{c}$ and the intersections of ${\mathcal C}$ with maximal inexact subsets of ${\mathcal A}$.
We say that ${\mathcal X}\subset {\mathcal A}_{c}$ is a {\it complementary} subset if 
${\mathcal A}_{c}\setminus {\mathcal X}$ is a maximal inexact subset of ${\mathcal A}_c$.
Complementary subsets of ${\mathcal A}_{c}$ can be characterized as the intersections of ${\mathcal C}$ with complementary subsets of ${\mathcal A}$.
 
 A bijection $g$ from an apartment of ${\mathcal C}$ to an apartment of ${\mathcal C}'$ is said to be {\it special} 
if $g$ and $g^{-1}$ transfer inexact subsets to inexact subsets.
For every special bijection one of the following possibilities is realized:
\begin{enumerate}
\item[(1)] The intersections of ${\mathcal C}$ with simple subsets of first and second type go to 
the intersections of ${\mathcal C}'$ with simple subsets of first and second type, respectively.
\item[(2)] The intersections of ${\mathcal C}$ with simple subsets of first type go to the intersections of ${\mathcal C}'$ with simple subsets of second type
and the intersections of ${\mathcal C}$ with simple subsets of second type correspond to the intersections of ${\mathcal C}'$ with simple subsets of first type.
In this case, we have ${\mathcal G}={\mathcal G}_{\alpha}(V)$.
\end{enumerate}
The proof is similar to the proof of Lemma \ref{lemma1-1}.

Then $f$ is adjacency preserving in both directions, 
i.e. $f$ is an isomorphism between the restrictions of $\Gamma({\mathcal G})$ to ${\mathcal C}$ and ${\mathcal C'}$
(the proof is similar to the proof of Lemma \ref{lemma1-3}).
As in the proof of Lemma \ref{lemma1-4}, we show that $f$ sends stars to stars and tops to tops.
Then, by Subsection 3.1, $f$ can be uniquely extended to a bijection ${\overline f}:{\mathcal C}_{\pm}\to {\mathcal C}'_{\pm}$
preserving the inclusion relation in both directions.

If ${\mathcal A}$ is an apartment of ${\mathcal G}$ and ${\mathcal A}\cap {\mathcal C}$ is an apartment of ${\mathcal C}$,
then $f({\mathcal A}\cap {\mathcal C})$ is an apartment of ${\mathcal C}'$ and
$$f({\mathcal A}\cap {\mathcal C})={\mathcal A}'\cap {\mathcal C}',$$
where ${\mathcal A}'$ is a certain apartment of ${\mathcal G}$.
Since ${\overline f}:{\mathcal C}_{\pm}\to {\mathcal C}'_{\pm}$ is inclusion preserving in both directions, we have 
$${\overline f}({\mathcal A}\cap {\mathcal C}_{\pm})={\mathcal A}'\cap {\mathcal C}'_{\pm}.$$
Therefore, ${\overline f}$ sends the intersections of ${\mathcal C}_{\pm}$ with apartments of ${\mathcal G}$ to 
the intersections of ${\mathcal C}'_{\pm}$ with apartments of ${\mathcal G}$.
The similar statement holds for the inverse bijection.

Using the latter property and the arguments from Subsection 4.4, 
we establish that for every $X\in {\mathcal C}_{\pm}$ there is a semilinear isomorphism 
$l_{X}:X\to {\overline f}(X)$ such that ${\overline f}(Y)=l_{X}(Y)$ for every $Y\in {\mathcal C}_{\pm}$ contained in $X$.
Let $h_{X}$ be the isomorphism of $\Pi_{X}$ to $\Pi_{{\overline f}(X)}$ induced by $l_{X}$.
We prove the direct analogue of Lemma \ref{lemma2-1} which implies the existence of an automorphism $h$ of $\Pi_{V}$
such that the restriction of $h$ to every ${\mathcal G}_{1}(X)$, $X\in {\mathcal C}_{\pm}$ coincides with $h_{X}$.
This completes the proof.

\end{document}